\documentstyle[11pt]{article}
\def\vs{\vspace*{3mm}}

\def\ms{\; \; \;}
\def\rno{\mbox{I\hspace*{-.7mm}R}}
\def\be{\begin{equation}}
\def\ee{\end{equation}}
\def\bc{\begin{center}}
\def\ec{\end{center}}
\def\bfl{\begin{flushleft}}
\def\efl{\end{flushleft}}
\def\bfr{\begin{flushright}}
\def\efr{\end{flushright}}
\def\bi{\begin{itemize}}
\def\ei{\end{itemize}}
\def\ben{\begin{enumerate}}
\def\een{\end{enumerate}}
\def\bt{\begin{tabular}}
\def\et{\end{tabular}}
\def\d#1{{#1\kern-0.4em\char"16\kern-0.1em}}
\def\D#1{{\raise0.2ex\hbox{-}\kern-0.4em #1}}
\def\zn{,\kern-0.1em ,}
\def\eps{\varepsilon}
\newtheorem{th}{Theorem}
\newtheorem{lem}{Lemma}
\title{\Large {\bf PIECEWISE EQUIDISTANT MESHES FOR QUASILINEAR TURNING 
POINT PROBLEMS: TECHNICAL REPORT}\footnote{This manuscript was prepared in February, 2000, but has
remained unpublished and I decided to make it accessible through arXiv.} 
}
\author{{\sc Relja Vulanovi\'c}\\
{\scriptsize Kent State University at Stark}\\
{\scriptsize 6000 Frank Ave NW, North Canton, OH 44720, USA}}
\date{}
\begin{document}
\maketitle
\begin{abstract}
A class of quasilinear singularly perturbed boundary value problems with a 
turning point of attractive type is considered. The problems are solved 
numerically by a finite--difference scheme on a special discretization mesh 
which is dense near the turning point. The scheme is a combination of the
standard central and midpoint schemes and is practically second--order accurate.
Pointwise accuracy is uniform in the perturbation parameter, $\eps$, and, 
moreover, $L^1$ errors decrease when $\eps\to 0$. This is achieved by the use 
of meshes which generalize the piecewise equidistant Shishkin mesh. Two 
particular types of meshes are considered and compared. 
\vs
\end{abstract}

\section{Introduction}
The following singularly perturbed boundary value problem was considered in
Vulanovi\'c and Lin \cite{rvl}:
%1
\be -\eps u'' -xb(x,u)u'+c(x,u) =0, \ms x\in [-1,1] , \ms u(\pm 1)=U_{\pm} , \ee
where $\eps$ is a small positive parameter, $b>0$ and $c=O(|x|+\eps)$ are 
sufficiently smooth functions, and $U_\pm$ are two given numbers. 
Some additional conditions were assumed as well, but we 
shall not list them here. It was shown that any solution of problem (1) has  
an interior layer of exponential type at the turning
point $x=0$. An appropriate numerical method was proposed, based on finite 
differences on a special discretization mesh dense in the layer. The error of 
this method was estimated in a discrete $L^1$ norm by 
%2
\be M(\sqrt{\eps}+e^{-N})N^{-1} . \ee
(In (2) and throughout the present paper, $N$ denotes the number of mesh steps 
and $M$ stands for any, in the sense of $O(1)$, positive constant independent 
of $\eps$ and $N$.)

In the present paper, we are interested in improving the result from \cite{rvl}.
This is achieved by applying a different kind of discretization mesh and 
a higher--order scheme. To simplify the presentation, we shall only consider a
special case of problem (1), viz.\ the case $b=b(u)$ and $c\equiv 0$. However, 
the same theoretical results hold for the general problem if the method and
some conditions are modified appropriately.

Singularly perturbed boundary value problems arise in various applications,
see Chang and Howes \cite{chh} for instance. Two recent books discuss numerical 
methods for these problems: Miller et al.\ \cite{mil} and Roos et al.\ 
\cite{rst}. Numerical methods for different types of turning point problems have
attracted a considerable attention, let us only mention the papers by Berger 
et al.\ \cite{berg}, Lin \cite{lin}, Vulanovi\'c and Farrell \cite{rvf}, 
Clavero and Lisbona \cite{cl}, and Sun and Stynes \cite{ss}. Prior to 
\cite{rvl}, some weaker versions of problem (1) were considered in 
Vulanovi\'c \cite{rvbail}, \cite{rvman}, and \cite{rvim}, where special 
discretization meshes were also used. The meshes belong to the class of 
explicitly constructed meshes, which means that they are formed before the 
discrete problem is solved. This requires sharp derivative estimates of the 
continuous solution. The same approach will be applied here.

In general, the explicitly constructed meshes can be divided in two main 
classes, meshes of Bakhvalov type (B meshes) and meshes of Shishkin type
(S meshes). The former were introduced in Bakhvalov \cite{bakh} and later on
generalized and simplified in Vulanovi\'c \cite{rvzb}. Many different types
of singular perturbation problems have been successfully solved numerically on 
B meshes. A B mesh is constructed by a smooth mesh generating function
which maps equidistant points into mesh points which are dense in the layer(s).
In each layer, the mesh generating function corresponds to the inverse of
the function describing how the solution behaves there. On the other hand,
S meshes (Shishkin \cite{sh}) are piecewise equidistant, and therefore
much simpler. An S mesh for a problem like (1) would typically consist
of three  equidistant parts: a fine part around $x=0$ and two coarse parts of
the mesh outside the layer, with the transition points between these parts 
located at $\pm\alpha\sqrt{\eps}\ln N$, where $\alpha>0$ is independent of 
$\eps$ and of $N$.
However, S meshes produce less accurate numerical results (see
Vulanovi\'c \cite{rvh} for a comparison of B and S meshes). There are two
possibilities for improving numerical results obtained on S meshes. One of them 
is to use more accurate schemes that may be more complicated but still easier to
analyze on S meshes than on B meshes, and another one is to improve the mesh
itself. This paper uses both approaches: a scheme which is practically
second--order accurate and a modified and improved S mesh.

Different improvements of S meshes have been considered so far. That by
Lin{\ss} \cite{tor1} uses a modification which makes the mesh more similar to
a B mesh. The resulting mesh is not piecewise equidistant any more. The
approach in Vulanovi\'c \cite{rvpre} improves 
the S mesh while keeping it piecewise equidistant (the mesh has more
equidistant parts which are constructed in a special way). A recent paper
by Roos and Lin{\ss} \cite{rl} (see Lin{\ss} et al.\ \cite{lrv} as well) 
provides for a unified theory which covers both S and B meshes (the latter only 
slightly different from those in \cite{rvzb}) and their generalizations. All 
these papers deal with non--turning point problems. 

We are not going to consider all the different types of meshes here. We shall 
only analyze the mesh from \cite{rvpre} and its slight modification. 
In this way, we show that meshes of this type can be applied also to 
turning point problems. The modification of the mesh from \cite{rvpre} 
consists of replacing $N$ in the transition
point formulas by $1/\eps$. Such transition points in S meshes were briefly 
discussed in \cite{rl} in the non--turning point case. They are closer to the
points marking the beginning and the end of the layer and they 
improve both theoretical and numerical results for problems of type (1).
The two kinds of meshes give quite 
satisfactory results and the discretization scheme is easier to discuss on 
these piecewise equidistant meshes. This is why we are not going to consider
B meshes here. Besides, the particular B mesh applied to (1) in \cite{rvl} is 
even more complicated than some other, more standard, B meshes.

The scheme which we shall use in this paper belongs to the class of hybrid
(or switching) finite--difference schemes. The standard central scheme is used 
inside the layer, where the mesh is fine and the scheme is unconditionally 
stable, whereas a midpoint upwind scheme is used outside the layer, on the 
coarse part of the mesh. For such schemes, see Vulanovi\'c \cite{rvgs} and
\cite{rveo}, and more recently, Stynes and Ross \cite{sr}, and Lin{\ss}
\cite{tor2}. Of these papers, only \cite{rveo} and \cite{tor2} deal with
hybrid schemes for quasilinear problems. The scheme in \cite{tor2} is less 
general, since it is constructed for a non--turning point problem. However, it 
is simpler and that is why we are going to use a very similar approach here. 

The paper is organized as follows. Precise assumptions on the continuous
problem and properties of its solution are given in section 2. In addition to 
the interior turning point case, a boundary turning point is also considered.
In that case, the turning point is still $x=0$, but the interval $[0,1]$ is
considered instead of $[-1,1]$. The numerical method is easier to describe 
for the boundary turning point problem, thus this case also serves the 
purpose of simplifying the presentation in section 3. Subsection 3.1
introduces the discretization scheme and analyzes its stability.
The special mesh is described and the main result is stated and 
proved in subsection 3.2. The necessary changes for the interior turning point 
case are explained in section 4. Some additional remarks are also given there.
This is followed by numerical results in section 5. 

To illustrate our main result, let us state it for the interior
turning point case when the mesh consists of five equidistant parts (the central
one around $x=0$ being the finest). Then, an error estimate of the form 
\[ M\sqrt{\eps}\left[\ln\left(\ln\frac{1}{\eps}\right)\right]^2 N^{-2} \]
can be proved in a discrete $L^1$ norm. This is an improvement over (2).

Let us finally mention that \cite{ss} is to our knowledge the only other paper 
which uses a piecewise equidistant mesh to solve a turning point problem
numerically. However, the problem considered there is different from (1)
and requires a different, more complicated mesh. In particular, the number
of equidistant parts of the mesh depends on $N$, which is not the case here.

\section{The Continuous Problem}
For simplicity, we are going to use $\eps^2$ instead of $\eps$ in the rest of 
the paper. We consider the problem
%3
\be \eps^2 u'' +  a(x,u)u' =0, \ms x\in I=[\nu,1], \ms u(\nu)=U_-, 
\ms u(1) =U_+ , \ee
where $\eps$ is a perturbation parameter, $0<\eps<<1$, and where $\nu$ stands 
for either 0 or $-1$. We assume that
%4
\be a(x,u)=xb(u) \ee
with $b\in C(\rno)$, 
and that $U_\pm$ are two different constants (otherwise $U_+=U_-$ solves (3)). 
The case $\nu=0$ describes a boundary turning point problem, whereas if 
$\nu=-1$, we have an interior turning point problem.

We can assume without loss of generality that $U_-<U_+$. Then $U_-$ and
$U_+$ are respectively the lower and upper solutions of (3),
and therefore, the problem has a unique solution, $u_\eps\in C^2(I)$, satisfying
\[  u_\eps(x)\in U:=[U_-, U_+] , \ms x\in I \]
(see Lorenz \cite{lor}).
Note that $U_\pm$ are also the solutions of the reduced equation
\[ a(x,u)u'=0 , \ms x\in I, \]
subject to only one of the original boundary conditions.

Throughout the paper, we shall assume that 
$b\in C^3(U)$, so that $u_\eps\in C^5(I)$.
Another assumption that will be needed in this paper is
%5
\be b^*\ge b(u)\ge b_* > 0 , \ms u\in U,  \ee
(of course, the upper bound on $b$ is not a restriction here).
Then the  solution $u_\eps$ of the problem (3)--(5) satisfies the following 
estimate:
%6
\be |u_\eps^{(k)}(x)|\le M\eps^{-k} y(x)\le M\eps^{-k} z(x), 
\ms k=0,1,2,3,4, \ms x\in I , \ee
where
\[ y(x)= e^{-b_{**}x^2/2\eps^2}  \]
and
\[ z(x) = e^{-m|x|/\eps} . \]
The above constants $b_{**}$ and $m$ are positive and independent of $\eps$. 
$b_{**}$ satisfies $b_{**}<b_*$, whereas $m$ is arbitrary.
As the estimate (6) is sharp, it clearly shows that $u_\eps$ has a layer
of exponential type at $x=0$. Moreover,
%7
\be |u_\eps(x)-U_+| \le M y(x)\le M z(x), \ms x\in [0,1] , \ee
and similarly
\[ |u_\eps(x)-U_-| \le M y(x) \le M z(x), \ms x\in [-1,0] , \ms 
\mbox{if $\nu=-1$} . \]
How to prove (6) and (7), can be found
in \cite{rvl} and \cite{rvbail}, cf.\ \cite{rvman} as well. Note that the proof 
of (6) requires $b\in C^3(U)$.  

It is another novelty of this paper, as compared to \cite{rvl}, 
\cite{rvbail}, and \cite{rvim}, that the condition (5) is given 
locally, i.e.\ for $u\in U$, and not for $u\in\rno$. 

\section{The Case $\nu=0$}
\subsection{The Discretization}
Let $I^h$ denote the discretization mesh with points $x_i$, $i=0,1,\ldots,N$,
$0=x_0<x_1<\cdots<x_N=1$, and let $h_i=x_i-x_{i-1}$, $i=1,2,\ldots,N$. The
only assumption on $I^h$ needed here is 
%8
\be h_i\le h_{i+1}, \ms i=1,2,\ldots, N-1 , \ee
(the special mesh will be introduced in the next subsection).
Also, let $\hbar_i=(h_i+h_{i+1})/2$, 
$i=1,2,\ldots,N-1$ and $x_{i\pm\frac 12}=(x_i+x_{i\pm 1})/2$.
Let $w^h$ be an arbitrary mesh  function on
$I^h\setminus \{0,1\}$, which is identified with a column vector in 
$\rno^{N-1}$,
\[ w^h=[w_1,w_2,\ldots,w_{N-1}]^T , \]
where for simplicity $w_i=w_i^h$. For any mesh function, we shall formally
set $w_0=U_-$ and $w_N=U_+$. The restriction of the continuous solution
$u_\eps$ on $I^h\setminus\{0,1\}$ will be denoted by $u_\eps^h$. 
Let $W=\{ w^h\in \rno^{N-1} \; |\; w_i\in U, \ms i=1,2,\ldots,N-1\}$.

The following quantity is used to define the discretization:
\[ \rho_i = \frac{b^*x_{i-1}h_i}{2\eps^2} , \]
where $b^*$ is given in (5). Let us consider the set of indices
\[ J=\{ i \; | \; \rho_i\le 1 \} \subseteq \{1,2,\ldots,N-1\} . \]
Note that $1\in J$ and let $n=\max J$. Because of (8),
%9
\be \rho_i \le 1 , \ms i=1,2,\ldots,n . \ee
If $1<n<N-1$, we define
\[ \chi_i = \left\{\begin{array}{ll}
	\hbar_i & \mbox{if}\ms 1\le i \le n-1 , \\
	\frac{\textstyle h_i}{\textstyle 2}+h_{i+1} & \mbox{if}\ms i=n , \\
	h_{i+1} & \mbox{if}\ms n+1\le i\le N-1 .
\end{array}\right. \]
If $n=1$, $\chi_i=h_{i+1}$ and if $n=N-1$, $\chi_i =\hbar_i$, in both cases 
for all $i=1,2,\ldots, N-1$. 

We can now introduce the finite--difference operators 
\[ D''w_i = \frac{1}{\chi_i}\left(\frac{w_{i-1}-w_i}{h_i}+
\frac{w_{i+1}-w_i}{h_{i+1}}\right), \]
\[ D'w_i = \frac{w_{i+1}-w_{i-1}}{2\hbar_i} , \]
\[ D'_+ w_i = \frac{w_{i+1}-w_i}{h_{i+1}} , \]
\[ D'_{t} w_i = \frac{2w_{i+1}-w_i-w_{i-1}}{2\chi_i} , \]
\[ D^\circ w_i = \frac{w_{i+1} + w_i}{2} . \]
The differential equation of problem (3) is discretized in the following form:
%10
\be -\eps^2 u'' - f(x,u)'+ f_x(x,u) = 0 , \ee
where 
\[ f(x,u)=\int_{U_+}^u a(x,t)dt . \]
When discretizing $f$ on $I^h$, we use the notation $f_i=f(x_i,w_i)$. 
The notation $f_{x,i}$, $a_i$, $b_i$, etc.\ has an analogous meaning. Then 
the following schemes are used to discretize (10):
\[ T_c w_i= -\eps^2 D''w_i - D'f_i + f_{x,i} , \]
\[ T_+ w_i = -\eps^2 D'' w_i - D'_+ f_i +D^\circ f_{x,i} , \]
\[ T_t w_i = -\eps^2 D''w_i - D'_t f_i + f_{x,i} . \]
By combining those schemes, we obtain the discretization of (10) on $I^h$,
%11
\be T w^h = 0 , \ee
where, if $1<n<N-1$,
\[ Tw_i = \left\{ \begin{array}{ll}
	T_c w_i & \mbox{if}\ms 1\le i\le n-1 ,\\
	T_t w_i & \mbox{if}\ms i=n ,\\
	T_+ w_i & \mbox{if}\ms n+1\le i \le N-1 .
\end{array}\right. \]
$T_c$ is the standard central scheme, $T_+$ is the midpoint scheme used
to discretize (10) at the point $x_{i+\frac 12}$, and $T_t$ is
a transition scheme between $T_c$ and $T_+$. If $n=1$, $T_c$ is not used
and there is no need for the transition scheme. In that case, $T\equiv T_+$. 
Likewise, if $n=N-1$, $T_+$ is not used and we set $T\equiv T_c$.
For simplicity, in what follows, we shall only consider $1<n<N-1$.

The central and midpoint schemes
are combined above in the same way as in \cite{rvgs}, but the transition scheme
was not required for the type of problems considered there. The present 
transition scheme is a little simpler than the one used in \cite{tor2}. Our 
stability analysis needs such a transition, but this may be just a technical 
requirement. Note that $\chi_i=h_{i+1}$ when $T_+$ is used. This gives a 
nonstandard scheme $D''$ for discretizing $u''$. However, such schemes have been
used earlier, see \cite{rvim} and \cite{lrv}.

Let us finally introduce some vector and matrix norms.
By $\|\cdot\|_\infty$ and $\|\cdot\|_1$ we denote the vector norms
\[ \|w^h\|_\infty = \max_{1\le i\le N-1}|w^h|, \ms 
\|w^h\|_1 = \sum_{i=1}^{N-1} |w_i| , \]
and, at the same time, their subordinate matrix norms. The diagonal matrix
\[ H=\mbox{diag}(\chi_1,\chi_2,\ldots,\chi_{N-1}) \]
is used to define the following discrete $L^1$ norm:
\[ \|w^h\|_H = \|Hw^h\|_1 , \]
and its subordinate matrix norm
\[ \|A\|_H = \|HAH^{-1}\|_1 , \]
where $A$ is an arbitrary $(N-1)\times(N-1)$ matrix. The vector norm
\[ \sum_{i=1}^{N-1}\hbar_i|w_i| \]
and the corresponding matrix norm are usually used for nonequidistant
discretizations of quasilinear problems, see \cite{rvl} for instance.
For the modified norms like $\|\cdot\|_H$ above, cf.\
\cite{rvim} and \cite{tor2}.

Let $G=[g_{ij}]=T'(w^h)$ be the 
Fr\'echet derivative of the discrete operator $T$ on mesh $I^h$ at some 
$w^h\in W^h$. 

\begin{lem} Let (5) and (8) hold. Then $G$ is an $L$--matrix.
\end{lem}
{\sc Proof.} Since $G$ is a tridiagonal matrix, we have to show that
\[ g_{ii} > 0 \ms \mbox{and}\ms g_{i,i\pm 1} \le 0 . \]
It is easy to see that $g_{ii}>0$ for all the schemes used. When $T_c$ or
$T_t$ are applied, (9) guarantees that $g_{i,i\pm 1}\le 0$.
For $T_+$, $g_{i,i-1}\le 0$ is immediate, and
\[ g_{i,i+1}  =  
-\frac{\eps^2}{h_{i+1}\hbar_i} -\frac{a_{i+1}}{h_{i+1}}+ \frac{b_{i+1}}{2} 
 \le  b_{i+1}\left(\frac 12-\frac{x_{i+1}}{h_{i+1}}\right) < 0 , \]
since $h_{i+1} < x_{i+1}$. \hfill $\Box$

\begin{th} Let (5) and (8) hold. Then the discrete problem (11) has a unique 
solution, $w^h_\eps$, which belongs to $W$. Moreover,
the following stability inequality holds for any two mesh functions
$w^h$, $v^h\in W$:
%12
\be \|w^h-v^h\|_H \le \frac{2}{b_*}\|Tw^h-Tv^h\|_H . \ee
\end{th}
{\sc Proof.} Let $e^h=[1,1,\ldots,1]^T$ $\in \rno^{N-1}$. It can be shown that
%13
\be s^h:= \left( HGH^{-1}\right)^Te^h \ge \frac{b_*}{2}e^h , \ee
where the inequality should be understood componentwise. The proof of (13) is 
elementary and it requires the transition scheme. In fact, if $i\ne n+1$,
\[ s_i = \frac{\chi_{i-1}}{\chi_i}g_{i-1,i}+ g_{ii} + 
\frac{\chi_{i+1}}{\chi_i}g_{i+1,i}\ge b_i \ge b_* , \]
where we formally set $g_{01}=g_{N,N-1}=0$. If $i=n+1$, then from
$T_t$ we get
\[ \frac{\chi_{i-1}}{\chi_i}g_{i-1,i} = -\frac{\eps^2}{h_i h_{i+1}}-
\frac{a_i}{h_{i+1}} , \]
and $T_+$ gives
\[ g_{ii} = \eps^2\frac{2\hbar_i}{h_i h^2_{i+1}}+\frac{a_i}{h_{i+1}}+
\frac{b_i}{2} , \]
and
\[ \frac{\chi_{i+1}}{\chi_i}g_{i+1,i} = -\frac{\eps^2}{h^2_{i+1}}, \]
so that
\[ s_i = \frac{b_i}{2} \ge \frac{b_*}{2} . \]
The discussion above also illustrates how the transition scheme enables
the proof of (13).

The inequality (13) implies that $G$ is an inverse--monotone matrix (and therefore an 
M--matrix) and also that
\[ \| G^{-1}\|_H \le \frac{b_*}{2} . \]
This result can be applied immediately to the matrix 
\[ A=\int_0^1 T'(v^h+s(w^h-v^h)) ds \] 
in
\[ Tw^h-Tv^h = A(w^h-v^h) , \]
and (12) follows.        

That (11) has a solution in $W$ can be proved by showing that
\[ T (U_- e^h) \le 0 \le T (U_+ e^h) . \]
The second inequality is immediate because of the way $f$ is defined.
To illustrate the proof of the first inequality, let us consider
\[ [T(U_-e^h)]_n = [T_t(U_-e^h)]_n = \int_{U_+}^{U_-}b(t)\left(1-
\frac{2x_{n+1}-x_n-x_{n-1}}{h_n+2h_{n+1}}\right)dt = 0 , \]
and
\begin{eqnarray*}
 [T(U_-e^h)]_{N-1} & = & [T_+(U_-e^h)]_{N-1} \\ 
 & = & -\eps^2\frac{U_+-U_-}{h_N^2}+
\int_{U_+}^{U_-} b(t)\left(\frac{x_{N-1}}{h_N}+\frac 12\right) < 0 
\end{eqnarray*}
(recall that $w_N$ is replaced by $U_+$). 

The solution is unique because of (12). \hfill $\Box$

\subsection{The Main Result}
We shall now define the special discretization mesh. Let $\lambda$ be either
$1/\eps$ or $N$ and let
\[ \ln^0 \lambda = \lambda , \ms
\ln ^k \lambda = \ln(\ln^{k-1}\lambda), \ms k=1,2,\ldots, K, \]
where $K=K(\lambda)$ is a positive integer such that $0<\ln^K\lambda<1$. Let
$\ell$, $1\le \ell\le K$, denote a fixed integer independent of $\eps$.
Also, let $\alpha$ be a positive constant independent of $\eps$ and $N$.
Then we define the transition points
\[ \tau_k=\alpha\eps\ln^{\ell-k+1}\lambda , \ms k=1,2,\ldots, \ell . \]
We shall assume that $\tau_\ell<1$, since $\eps<<1$. By formally setting 
$\tau_0=0$ and $\tau_{\ell+1}=1$, we can split up the interval $[0,1]$ into 
$\ell+1$ subintervals,
\[ [0,1]=\bigcup_{k=1}^{\ell+1} I_k, \ms I_k=[\tau_{k-1}, \tau_k] . \]  
Each interval $I_k$ is then divided into $N_k\ge 2$ equidistant subintervals, 
so that
\[ N_1+N_2+\ldots + N_{\ell+1}=N  \]
and
\[  N \le M N_k, \ms k=1,2,\ldots,\ell+1 . \]
The points obtained in this way form the mesh on $[0,1]$, which we shall 
refer to as the S($\ell$) mesh (Shishkin mesh with $\ell$ transition points).
The standard S mesh is S(1) with $\lambda=N$. S($\ell$) satisfies (8).

Note that if (9) holds with $n=N-1$ on the S($\ell$) mesh, this practically 
means that $1/N \le M\eps^2$, whereas usually $\eps\le 1/N$. Thus it is not
realistic to expect that the discrete operator $T$ be identical to the
central scheme $T_c$. On the other hand, it is possible that $n=1$ and 
$T\equiv T_+$ if $\eps$ is sufficiently small, 
$\eps^2 \le M N^{-2}(\ln^\ell\lambda)^2$.
Even though we are only showing details for $1<n<N-1$, all our results are 
true for $n=1$ or $n=N-1$ as well.

Let
\[ \delta = \eps\left(\frac{\ln^\ell\lambda}{N}\right)^2  \]
and
%14
\be v(x) = e^{-\beta x/\eps}  \ee
with a positive constant $\beta$.

We first consider the case $\lambda=1/\eps$.

\begin{lem} On the S($\ell$) mesh with $\lambda=1/\eps$ it holds that
%15
\be \frac{h^2_{i+1}}{\eps} v_{i-1} \le M \delta , \ee 
where $i=1,2,\ldots, N-1$ and $v_i=v(x_i)$ with a sufficiently large constant 
$\beta$ independent of $\eps$ and $N$.
\end{lem}
{\sc Proof.} We shall consider several cases.
If $x_i\in(0,\tau_1)$, then $h_{i+1}\le M\eps N^{-1} \ln^\ell\lambda$ and
(15) follows immediately. 

In all other cases (15) can be proved with 
$\eps N^{-2}$ instead of $\delta$.
If $x_i\in (\tau_k,\tau_{k+1})$ for $1\le k\le \ell$, then
\begin{eqnarray*} \frac{h^2_{i+1}}{\eps} v_{i-1} & \le & 
M \eps \left(\frac{\ln^{\ell-k}\lambda}{N}\right)^2 
e^{-\alpha\beta \ln^{l-k+1}\lambda}\\
& = & M\eps N^{-2} (\ln^{\ell-k}\lambda)^{2-\alpha\beta} \le M \eps N^{-2} ,
\end{eqnarray*}
where $\beta$ is chosen so that $\alpha\beta\ge 2$.

If $x_i=\tau_k$, $k=1,2,\ldots,\ell$, then 
\[ x_{i-1}\ge \tau_k\left(1-\frac{1}{N_k}\right)\ge \frac{\tau_k}{2} \]
and (15) follows like in the previous case but with $\alpha\beta\ge 4$.\hfill
$\Box$

Let us introduce some more notation. By $\cal I$ we denote
\[ {\cal I}=\bigcup_{k=0}^\ell {\cal I}_k , \]
where
\[ {\cal I}_k = \{ i \; |\; x_i\in (\tau_k,\tau_{k+1})\} . \]

\begin{lem} On the S($\ell$) mesh with $\lambda=1/\eps$ it holds that
\[ \sigma_k:=\sum_{i\in {\cal I}_k} \left(\frac{h_{i+1}}{\eps}\right)^2 
\int_{x_{i-1}}^{x_{i+1}} v(x) dx \le M\delta , \]
where $k=0,1,\ldots,\ell$ and where $v(x)$ is given in (14) with a sufficiently 
large constant $\beta$ independent of $\eps$ and $N$.
\end{lem}
{\sc Proof.} We have
\begin{eqnarray*} 
\sigma_k  & \le &  M\left(\frac{\tau_{k+1}}{\eps N}\right)^2
\int_{\tau_k}^{\tau_{k+1}} v(x) dx \\ 
& \le & M\left(\frac{\tau_{k+1}}{\eps N}\right)^2 \eps v(\tau_k) \le M \delta .
\end{eqnarray*}
The last inequality follows similarly to the proof of Lemma 2. \hfill $\Box$

We can now prove the main result of the paper.

\begin{th} 
Let (5) hold and let $w^h_\eps$ be the solution of the 
discrete problem (11) on the S($\ell$) mesh with $\lambda=1/\eps$. Then the 
following error estimate holds:
\[ \|w^h_\eps- u^h_\eps\|_H \le M \delta . \]
\end{th}
{\sc Proof.} Using (12) with $w^h=w^h_\eps$ and $v^h=u^h_\eps$, we see that
in order to prove the theorem, it suffices to show that the consistency error
$r^h=Tu^h_\eps$ satisfies 
%16
\be \|r^h\|_H \le M\delta. \ee
Let
\[ r^h_c=T_c u^h_\eps \ms \mbox{and} \ms r^h_+=T_+ u^h_\eps . \]
Since the transition between $T_c$ and $T_+$ can occur at any mesh point,
we can consider separately $\|r^h_c\|_H$ and $\|r^h_+\|_H$. As for $T_t$, it is
only used at $x_n$ and therefore, it is sufficient to estimate $\chi_n |r_n|$.

When the central scheme is used, Taylor's expansion of $r_i=r_{c,i}$ about 
$x_i$ gives
\begin{eqnarray*}
 \chi_i |r_i| & \le & M\left\{ \eps^2\left[ h_{i+1}(h_{i+1}-h_i)
|u_\eps^{(3)}(x_i)| + h_{i+1}^2\int_{x_{i-1}}^{x_{i+1}} |u_\eps^{(4)}(x)| dx 
\right]\right. \\
& + & \left. h_{i+1}(h_{i+1}-h_i)|p''(x_i)|+h_{i+1}^2\int_{x_{i-1}}^{x_{i+1}} 
|p^{(3)}(x)|dx \right\} ,
\end{eqnarray*}
where $p(x)=f(x,u_\eps(x))$ (for the integral form of the consistency error,
cf.\ Kellogg and Tsan \cite{kts} for instance). Then (6) implies
\[ \chi_i |r_i| \le M(Q_i+R_i) , \]
where
\[ Q_i=h_{i+1}(h_{i+1}-h_i)\eps^{-1} v_{i-1}  \]
and
\[ R_i=\left(\frac{h_{i+1}}{\eps}\right)^2\int_{x_{i-1}}^{x_{i+1}} v(x) dx . \]
Here, $v(x)$ is like in (14), with a constant $\beta$ satisfying $0<\beta<m$. 
$\beta$ is independent of $\eps$ and $N$, it can be chosen arbitrarily close to 
$m$, and therefore it can be made sufficiently large. On writing
\[ \|r^h\|_H = \sum_{i\not\in {\cal I}}\chi_i|r_i|+\sum_{i\in {\cal I}}
\chi_i|r_i| , \]
from Lemma 3 we immediately get
\[ \sum_{i\in {\cal I}}\chi_i|r_i|\le M \sum_{i\in {\cal I}} R_i \le M\delta . 
\]
On the other hand, if $i\not\in {\cal I}$, i.e.\ if $x_i=\tau_k$, it holds that
\[ Q_i \le \frac{h_{i+1}^2}{\eps} v_{i-1} \le M\delta , \]
because of Lemma 2. Since the number $\ell$ of transition points does not
depend on $N$, this proves (16) for $r^h=r^h_c$.

Let us now consider the consistency error when $T=T_+$. Then, $r_i=r_{+,i}$ is
expanded about $x_{i+1/2}$. The errors due to $D'_+$ and $D^\circ$ can be 
treated as above, but the one due to $D''$ requires a closer attention. Let
\[ r''_i = \eps^2[D''u_\eps(x_i) - u''_\eps(x_{i+1/2})]. \]
We have
\[ \chi_i|r''_i| \le M \eps^2\left[ (h_{i+1}-h_i)|u''_\eps(x_{i+1/2})| +
h_{i+1}\int_{x_{i-1}}^{x_{i+1}} |u^{(3)}_\eps(x)| dx \right] . \]
Since in this case $\rho_i>1$, it holds that
\[ \eps^2 \le Mx_{i-1}h_i . \]
This and (6) imply
\[ \chi_i|r''_i| \le M (Q''_i + R''_i) , \]
\[ Q''_i = h_i(h_{i+1}-h_i)\eps^{-1}v_{i-1} , \]
\[ R''_i = h^2_{i+1} \frac{x_{i-1}}{\eps^3}\int_{x_{i-1}}^{x_{i+1}}z(x) 
dx . \]
We can handle $Q''_i$ in the same way as $Q_i$ above, and
\begin{eqnarray*}
 R''_i & \le & \left(\frac{h_{i+1}}{\eps}\right)^2\frac{x_{i-1}}{\eps}
e^{-mx_{i-1}/2\eps}\int_{x_{i-1}}^{x_{i+1}} e^{-mx/2\eps} dx \\    
& \le & M \left(\frac{h_{i+1}}{\eps}\right)^2\int_{x_{i-1}}^{x_{i+1}}
e^{-mx/2\eps} dx ,
\end{eqnarray*}
so that Lemma 3 can be applied. Thus, (16) holds with $r^h=r^h_+$.    

Finally, let as consider the transition scheme at $x_n$. We again expand $r_n$ 
about $x_n$ to get
%17
\begin{eqnarray} 
\chi_n|r_n| & \le & M\left\{ \eps^2\left[ h_{n+1}|u''_\eps(x_n)| 
+h^2_{n+1}\max_{x_{i-1}\le x\le x_{i+1}}|u^{(3)}_\eps(x)|\right] \right. 
\nonumber \\
 & & \\
 & + & \left. h_{n+1}^2 |p''(\theta)| \right\} \nonumber
\end{eqnarray}
with some $\theta\in(x_{n-1},x_{n+1})$. The fact that in this case
$\rho_{n+1}>1$ implies
\[ \eps^2 \le Mx_nh_{n+1} . \]
We need this inequality and (6) in the estimate
\[ \eps^2 h_{n+1}|u''_\eps(x_n)| \le M h^2_{n+1}\frac{x_n}{\eps^2}z(x_n)
\le M \frac{h^2_{n+1}}{\eps} v_n . \]
The other terms on the right--hand side of (17) can be estimated directly
using (6). Thus, it follows that
\[ \chi_n|r_n| \le M \frac{h^2_{n+1}}{\eps}v_{n-1} \le M\delta , \]
where the last inequality follows from Lemma 2.\hfill $\Box$
\vs

Let us now turn to the case $\lambda=N$.

\begin{th}
Let (5) hold and let $w^h_\eps$ be the solution of the 
discrete problem (11) on the S($\ell$) mesh with $\lambda=N$. Then the 
following error estimate holds:
\[ \|w^h_\eps- u^h_\eps\|_H \le M \left(\delta + N^{-(1+\eta\ln N)}\right), \]
where $\eta$ is some positive constant independent of $\eps$ and $N$. 
\end{th}
{\sc Proof.} The result can be proved analogously to the proof of 
Theorem 2, cf.\ \cite{rvpre} as well. The only case
which produces the $N^{-(1+\eta\ln N)}$ term is when 
$x_i=\tau_\ell, \tau_\ell+\kappa$, 
where $\kappa$ is the mesh step in $[\tau_\ell,1]$. For instance, when 
estimating expressions like the left--hand side of (15) (which is also needed
in the proof of Lemma 3), the best we can get
at those points is 
\[ \frac{h_{i+1}^2}{\eps}v_{i-1} \le M \frac{1}{N^2\eps}e^{-\alpha\beta \ln N}
= M \frac{1}{\eps N^{2+\alpha\beta}} . \]
The case $x_i\in(\tau_\ell+\kappa,1)$ still gives the same result as before,
but requires a somewhat different technique:
\[ \frac{h_{i+1}^2}{\eps}v_{i-1} \le M \frac{\eps}{N^2}\frac{1}{N^{\alpha\beta}
\eps^2}e^{-\gamma/N\eps} \le M\frac{\eps}{N^2} , \]
where $\gamma$ is a positive constant independent of $\eps$ and $N$, and where
$\alpha\beta\ge 2$. Therefore, when $x_i=\tau_\ell, \tau_\ell+\kappa$, we use a
different estimate of the consistency error:
\[ \chi_i |r_i| = \chi_i |Tu_\eps(x_i)| \le M h_{i+1}\{ \eps^2
\max_{x_{i-1}\le x\le x_{i+1}} [|u_\eps''(x)| + |p'(x)| + |q(x)|]\} , \]
where $q(x)=f_x(x,u_\eps(x))$. We use here estimates (6) and (7) with
$y(x)$ to get
\[ \chi_i |r_i| \le M h_{i+1}y(x_{i-1})  \]
and 
\[ \chi_i |r_i| \le M N^{-1} e^{-\eta(\ln N)^2} \le M N^{-(1+\eta\ln N)} , \]
on setting $\eta=b_{**}\alpha^2/2$. \hfill $\Box$

The error estimate of Theorem 2 is better than that of Theorem 3
when $\eps\to 0$. It may look like the term $N^{-(1+\eta\ln N)}$ in the
estimate above arises for purely technical reasons. However, the numerical 
experiments in section 5 will show that when that term dominates in the error,
the error does not decrease together with $\eps$.  

\section{The Case $\nu=-1$ and Other Remarks}
Theorem 2 holds also for the case $\nu=-1$, when the problem (3) is considered
on $I=[-1,1]$, if the numerical method is modified appropriately. As for the
mesh, the simplest thing to do is to extend S($\ell$) to $[-1,0]$ symmetrically
to $x_0=0$: 
\[ x_{-i}=-x_i, \ms i=1,2,\ldots, N, \] 
\[ h_{-i}=x_{-i}-x_{-i-1}=h_{i+1} , \ms i=0,1,\ldots, N-1 . \]
This is accompanied with other
symmetrically changed definitions, like
\[ \chi_i = \left\{\begin{array}{ll}
	h_i & \mbox{if}\ms -N+1 \le i\le -n-1 , \\
	h_i+\frac{\textstyle h_{i+1}}{\textstyle 2} & \mbox{if}\ms i=-n, \\
	\hbar_i & \mbox{if}\ms -n+1\le i \le n-1 , \\
	\frac{\textstyle h_i}{\textstyle 2}+h_{i+1} & \mbox{if}\ms i=n , \\
	h_{i+1} & \mbox{if}\ms n+1\le i\le N-1 ,
\end{array}\right. \] 
and
\[ D'_- w_i = \frac{w_i-w_{i-1}}{h_i} , \]
\[ D'_{t-}w_i = \frac{w_{i+1}+w_i-2w_{i-1}}{2\chi_i} , \]
\[ D^\circ_- w_i = \frac{w_i+w_{i-1}}{2} . \]
The scheme is also symmetrical:
\[ Tw_i = \left\{ \begin{array}{ll}
	T_- w_i & \mbox{if}\ms -N+1\le i \le -n-1 \\
	T_{t-}w_i & \mbox{if}\ms i=-n, \\ 
	T_c w_i & \mbox{if}\ms -n+1\le i\le n-1 ,\\
	T_t w_i & \mbox{if}\ms i=n ,\\
	T_+ w_i & \mbox{if}\ms n+1\le i \le N-1 ,
\end{array}\right. \] 
where
\[ T_-w_i = -\eps^2 D''w_i -D'_- f_i + D^\circ_- f_{x,i} , \]
and
\[ T_{t-}w_i =-\eps^2 D''w_i - D'_{t-}f_i + f_{x,i} . \]
Note that $T_c$ is always used at least at $x_0=0$. To make the discretization
even more symmetric, we
also replace $U_+$ in $f$ with $U_-$ when discretizing (10) at $x_i\in(-1,0)$, 
and with $(U_+ + U_-)/2$ at $x_0=0$.

Under the conditions of Theorem 2, the following estimate holds:
%18
\be \| w^h_\eps -u^h_\eps \|_\infty \le M \frac{\ln^\ell (1/\eps)}{N} . \ee
This is because the smallest mesh step is $h_1=\tau_1/N$. The above result does
not mean $\eps$--uniform pointwise convergence, but $\ln^\ell(1/\eps)$ grows
very slowly when $\eps\to 0$. For instance, if $\ell=3$ and $\eps=10^{-12}$,
$\ln^\ell(1/\eps)\approx 1.2$. The numerical results of the next section show
that pointwise accuracy is even better than what (18) indicates. A similar
discussion is true for S($\ell$) with $\lambda=N$.

It is possible to adjust the present method to the most general problem (1)
under the conditions given in \cite{rvl}. The reduced solutions are not
$U_-$ and $U_+$ in that case, but some more complicated functions. They  
have to be incorporated in the function $f$, like it was done in \cite{rvl}.
In order to prove the stability inequality corresponding to (12),
the conditions of type (5) on $b$ have to hold for $u\in \rno$. The same
has to be assumed of all conditions on $c$. This is because in this case it is
generally difficult to find the upper and lower solutions of (1) and of the
discrete problem. For Theorems 2 and (3), more complicated
estimates corresponding to (6) and (7) have to be used, see \cite{rvl}.

\section{Numerical Results}
Our test problem is more general than (3) but less general than (1),
\[ -\eps^2 u'' - xuu'+ c(x)=0 , \ms u(\pm 1)= U_{\pm} , \]
where $c(x)$, $U_-\approx 1$, and $U_+\approx 3$ are determined by the exact 
solution being
\[ u_\eps (x)=2+\tanh\frac{x}{\eps} . \]

We have tested the S($\ell$) mesh with both $\lambda=1/\eps$ and $\lambda=N$,
and with values of $\ell=1,2,3$. Let
\[ q_k = \frac{N_k}{N}, \ms k=1,2,\ldots, \ell+1, \]
where $N$ is the number of mesh steps in $[0,1]$. 
Table 1 shows the values of the ratios $q_k$ that are used in the meshes below.
They were kept fixed for each $\ell$ regardless of other mesh parameters,
including $N$.
\vs
\bc
{\sc Table 1.} Mesh parameters 
\vs\\
\bt{|c||c|c|c|} \hline
$\ell$ & 1 & 2 & 3 \\ \hline\hline
$q_1$ & 3/4 & 1/4 & 1/8 \\ \hline
$q_2$ & -- & 1/2 & 1/8 \\ \hline
$q_3$ & -- & -- & 1/2 \\ \hline
\et
\ec
\vs
The choice of the ratios $q_k$ may influence the errors significantly. 
Other ratios have been also tested and the results for this test problem were 
the best when there were around  75\% of the mesh points in the layer.
The question of the optimal choice of the ratios is problem--dependent and
seems to be difficult to solve in general. 

The tables below show the errors
\[ E=E(N) =\|w^h_\eps -u^h_\eps\|_\infty \ms \mbox{and} \ms 
E_1=E_1(N) =\|w^h_\eps-u^h_\eps\|_H , \]
where the dependence on $N$ indicates that the mesh is used with $2N$ mesh
steps in $[-1,1]$. The corresponding numerical 
orders of convergence, Ord= Ord$(N)$ and Ord$_1$=Ord$_1(N)$, are also 
listed, where
\[ \mbox{Ord}(N) = \frac{\ln E(N) - \ln E(N/2)}{\ln 2} , \]
and Ord$_1(N)$ is defined analogously. 
The results in Tables 2--4 are given for 
the transition point coefficient $\alpha=1$, whereas Table 5 presents some
results for $\alpha=2$.

Table 2 illustrates that the errors are smaller if $\ell$ is larger,
that is, if the meshes have more subintervals within the layer. This is to be 
expected. The case $\lambda=1/\eps$ is shown. If $\lambda=N$, the errors behave 
analogously with respect to the change in $\ell$. 
\vs
\bc
\small
{\sc Table 2.} S($\ell$) mesh with $\lambda=1/\eps$, $\alpha=1$, $N=512$ 
\vs\\
\bt{|c||cc|cc|cc|cc} \cline{1-7}
$\ell$ & \multicolumn{2}{c|}{$\eps=2^{-14}$} & \multicolumn{2}{c|}
{$\eps=2^{-18}$} & \multicolumn{2}{c|}{$\eps=2^{-22}$} & & \\ \hline\hline
1 & 7.94--05 & 1.19--08 & 1.31--04 & 1.06--09 & 1.96--04 & 9.83--11 & 
$E$ & \multicolumn{1}{c|}{$E_1$} \\ 
  & 2.00 & 1.89 & 2.00 & 1.99 & 2.00 & 2.00 & Ord & \multicolumn{1}{c|}{Ord$_1$}
\\ \hline
2 & 3.84--05 & 6.91--09 & 4.74--05 & 4.15--10 & 5.53--05 & 3.00--11 & & \\ 
  & 2.00 & 1.81 & 2.00 & 1.98 & 2.00 & 2.00 & & \\ \cline{1-7}
3 & 1.96--05 & 4.29--09 & 2.29--05 & 1.93--10 & 2.60--05 & 1.33--11 & & \\ 
  & 2.00 & 1.68 & 2.00 & 1.96 & 2.00 & 2.00 & & \\ \cline{1-7}
\et
\ec
\vs

In the remaining tables, only $\ell=3$ is considered. Table 3 gives more
details of what can already be observed in Table 2, viz.\ errors $E_1$ 
decreasing together with $\eps$, while at the same time, errors $E$ slightly
increase, still preserving a high accuracy. The increase of $E$ is what we can 
expect from (18), but the accuracy and its order are higher than what (18)
indicates. Both Ord and Ord$_1$ are around 2, with Ord$_1$ being somewhat 
smaller for larger $\eps$ values.
\vs 
\bc
\small
{\sc Table 3.} S(3) mesh with $\lambda=1/\eps$, $\alpha=1$ 
\vs\\
\bt{|c||cc|cc|cc|cc} \cline{1-7}
$N$ & \multicolumn{2}{c|}{$\eps=2^{-14}$} & \multicolumn{2}{c|}
{$\eps=2^{-18}$} & \multicolumn{2}{c|}{$\eps=2^{-22}$} & & \\ \hline\hline
64 & 1.18--03 & 1.88--07 & 1.57--03 & 1.73--08 & 1.73--03 & 1.27--09 & 
$E$ & \multicolumn{1}{c|}{$E_1$} \\ 
  & 2.01 & 2.06 & 1.51 & 2.00 & 1.48 & 1.89 & Ord & \multicolumn{1}{c|}{Ord$_1$}
\\ \hline
128 & 3.13--04 & 4.80--08 & 3.63--04 & 2.97--09 & 3.94--04 & 2.45--10 & & \\ 
  & 1.92 & 1.97 & 2.12 & 2.54 & 2.13 & 2.37 & & \\ \cline{1-7}
256 & 7.83--05 & 1.38--08 & 9.16--05 & 7.51--10 & 1.04--04 & 5.33--11 & & \\ 
  & 2.00 & 1.80 & 1.99 & 1.98 & 1.92 & 2.20 & & \\ \cline{1-7}
512 & 1.96--05 & 4.29--09 & 2.29--05 & 1.93--10 & 2.60--05 & 1.33--11 & & \\ 
  & 2.00 & 1.68 & 2.00 & 1.96 & 2.00 & 2.00 & & \\ \cline{1-7}
\et
\ec
\vs

When $\lambda=N$, both $E$ and $E_1$ errors are more uniform in $\eps$ but
the accuracy is worse than for $\lambda=1/\eps$. This is shown in Table 4,
where we can see that Ord is less than 2 as can be expected from the
error estimate of Theorem 3. Also, $E_1$ is much worse than in 
Table 3 even though Table 4 shows Ord$_1$ significantly higher than 2. 
In fact, in Table 4,  $E_1$ does not
decrease together with $\eps$, which is possible according to the error
estimate of Theorem 3. This can be improved by increasing the value of
$\alpha$ so that the $\eps$--independent term of the error estimate gets
negligible. Table 5 shows that $E_1$ errors decrease for $\alpha=2$ when 
$\eps\to 0$ and that they become almost as accurate as in Table 3, but at the 
same time, the larger $\alpha$ spoils $E$ a little.
\vs 
\bc
{\sc Table 4.} S(3) mesh with $\lambda=N$, $\alpha=1$ 
\vs\\
\bt{|c||cc|cc|cc} \cline{1-5}
$N$ & \multicolumn{2}{c|}{$\eps=2^{-14}$} & \multicolumn{2}{c|}
{$\eps=2^{-18}$, $2^{-22}$} & & \\ \hline\hline
64 & 1.89--03 & 9.21--05 & 1.90--03 & 9.21--05 & 
$E$ & \multicolumn{1}{c|}{$E_1$} \\ 
  & 1.91 & 2.62 & 1.91 & 2.62 & Ord & \multicolumn{1}{c|}{Ord$_1$}
\\ \hline
128 & 5.07--04 & 1.41--05 & 5.10--04 & 1.42--05 & & \\ 
  & 1.90 & 2.70 & 1.89 & 2.70 & & \\ \cline{1-5}
256 & 1.37--04 & 2.10--06 & 1.39--04 & 2.10--06 & & \\ 
  & 1.89 & 2.76 & 1.88 & 2.75 & & \\ \cline{1-5}
512 & 3.64--05 & 3.02--07 & 3.76--05 & 3.04--07 & & \\ 
  & 1.91 & 2.79 & 1.88 & 2.79 & & \\ \cline{1-5}
\et
\ec
\vs
\bc
\small
{\sc Table 5.} S(3) mesh with $\lambda=N$, $\alpha=2$ 
\vs\\
\bt{|c||cc|cc|cc|cc} \cline{1-7}
$N$ & \multicolumn{2}{c|}{$\eps=2^{-14}$} & \multicolumn{2}{c|}
{$\eps=2^{-18}$} & \multicolumn{2}{c|}{$\eps=2^{-22}$} & & \\ \hline\hline
64 & 1.59--03 & 2.63--07 & 1.59--03 & 3.02--08 & 1.59--03 & 1.56--08 & 
$E$ & \multicolumn{1}{c|}{$E_1$} \\ 
  & 3.09 & 3.15 & .09 & 4.07 & 3.09 & 4.66 & Ord & \multicolumn{1}{c|}{Ord$_1$}
\\ \hline
128 & 4.44--04 & 4.81--08 & 4.44--04 & 3.51--09 & 4.44--04 & 7.24--10 & & \\ 
  & 1.84 & 2.45 & 1.84 & 3.10 & 1.84 & 4.43 & & \\ \cline{1-7}
256 & 1.40--04 & 1.67--08 & 1.40--04 & 1.07--09 & 1.40--04 & 8.44--11 & & \\ 
  & 1.66 & 1.53 & 1.66 & 1.73 & 1.66 & 3.10 & & \\ \cline{1-7}
512 & 4.34--05 & 5.53--09 & 4.34--05 & 3.46--10 & 4.34--05 & 2.23--11 & & \\ 
  & 1.69 & 1.59 & 1.69 & 1.61 & 1.69 & 1.92 & & \\ \cline{1-7}
\et
\ec
\vs

How fast do the $E_1$ errors in Tables 3 and 5 decrease as $\eps\to 0$?
This $\eps$--order can also be measured numerically, analogously to Ord$_1$:
\[ \mbox{Ord}_\eps = \frac{E_1(4\eps)- E_1(\eps)}{\ln 4} , \] 
where $N$ is kept fixed and the dependence of $E_1$ on the value of $\eps$ is
expressed. The results are given in Table 6. We can see that they confirm the
expected value of 1 for $\lambda=N$ and $\alpha=2$, whereas for 
$\lambda=1/\eps$, they are better than what Theorem 2 indicates, particularly 
for larger $\eps$.
\vs 
\bc
{\sc Table 6.} Ord$_\eps$ on S(3) mesh with $N=512$
\vs\\
\bt{|c||c|c|} \hline
$\eps$ & $\lambda=1/\eps$, $\alpha=1$ & $\lambda=N$, $\alpha=2$ \\ \hline\hline
$2^{-14}$ & 1.45 & 1.00 \\ \hline
$2^{-16}$ & 1.21 & 1.00 \\ \hline
$2^{-18}$ & 1.03 & 1.00 \\ \hline
$2^{-20}$ & 0.97 & 1.00 \\ \hline
$2^{-22}$ & 0.96 & 0.98 \\ \hline
\et
\ec

Let us finally mention that we have also tested the discretization in which
the central scheme is used instead of the transition schemes $T_t$ and $T_{t-}$.
The errors are somewhat worse, the difference being greater when 
$\lambda=1/\eps$ than when $\lambda=N$. Thus, it may be possible that the
use of the transition schemes is not entirely for technical reasons.


\begin{thebibliography}{10}

\bibitem{rvl}
R.~Vulanovi\'c and P.~Lin, Numerical solution of quasilinear attractive
turning point problems. {\em Computers Math.\ Applic.} {\bf 23}, 75--82
(1992).

\bibitem{chh}
K.~W.~Chang and F.~A.~Howes, {\em Nonlinear Singular Perturbation Phenomena:
Theory and Applications}, Springer, New York (1984).

\bibitem{mil}
J.~J.~H.~Miller, E.~O'Riordan, and G.~I.~Shishkin,
{\em Solution of Singularly Perturbed Problems with $\eps$-uniform
Numerical Methods -- Introduction to the Theory of Linear Problems in One
and Two Dimensions}, World Scientific, Singapore (1996).

\bibitem{rst}
H.-G.~Roos, M.~Stynes, and L.~Tobiska, {\em Numerical Methods for Singularly
Perturbed Differential Equations}, Springer, Berlin (1996).

\bibitem{berg} 
A.~E.~Berger, H.~Han, and R.~B.~Kellogg, A priori estimates and analysis
of a numerical method for a turning point problem. {\em Math.\ Comput.}
{\bf 42}, 465--492 (1984).

\bibitem{lin} 
P.~Lin, A numerical method for quasilinear singular perturbation problems with 
turning points. {\em Computing} {\bf 46}, 155--164 (1991).

\bibitem{rvf}
R.~Vulanovi\'c and P.~A.~Farrell, Continuous and numerical analysis of a
multiple boundary turning point problem. {\em SIAM J.\ Numer.\ Anal.}
{\bf 30}, 1400--1418 (1993).

\bibitem{cl}
C.~Clavero and F.~Lisbona, Uniformly convergent finite difference methods for
singularly perturbed problems with turning points. {\em Numer.\ Algorithms}\/
{\bf 4}, 339--359 (1993).

\bibitem{ss}
G.~Sun and M.~Stynes,
Finite element methods on piecewise equidistant meshes for interior
turning point problems. {\em Numer.\ Algorithms}\/ {\bf 8}, 111--129 (1994).  

\bibitem{rvbail}
R.~Vulanovi\'c, On numerical solution of some quasilinear turning point 
problems. In {\em Proc.\ BAIL V Conf.}, (Edited by Guo Ben--yu {\em et al.}),
pp.\ 368--373, Boole Press, Dublin (1988).

\bibitem{rvman}
R.~Vulanovi\'c, On numerical solution of a mildly nonlinear turning point
problem. {\em RAIRO Math.\ Model.\ Numer.\ Anal.} {\bf 24}, 765--784 (1990).

\bibitem{rvim}
R.~Vulanovi\'c, A uniform numerical method for a class of quasilinear
turning point problems. In {\em Proc.\ 13th IMACS World Congress on Computation 
and Applied Mathematics}, (Edited by J.J.H.~Miller and R.~Vichnevetsky),
p.\ 493, IMACS (1991).

\bibitem{bakh}
N.~S.~Bakhvalov,
Towards optimization of methods for solving boundary value problems
in the presence of a boundary layer. {\em Zh.\ Vychisl.\ Mat.\ Mat.\ Fiz.} 
{\bf 9}, 841--859 (1969), in Russian.

\bibitem{rvzb}
R.~Vulanovi\'c, On a numerical solution of a type of singularly perturbed
boundary value problem by using a special discretization mesh. {\em Univ.\ 
u Novom Sadu Zb.\ Rad.\ Prirod.--Mat.\ Fak.\ Ser.\ Mat.} {\bf 13}, 187--201
(1983).

\bibitem{sh}
G.~I.~Shishkin,
A difference scheme for a singularly perturbed parabolic equation
with a discontinuous boundary condition.
{\em Zh.\ Vychisl.\ Mat.\ Mat.\ Fiz.} {\bf 28}, 1679--1692 (1988), in Russian.

\bibitem{rvh}
R.~Vulanovi\'c, Forth order algorithms for a semilinear singular perturbation
problem. {\em Numer.\ Algorithms} {\bf 16}, 117--128 (1997).

\bibitem{tor1} T.~Lin{\ss},
An upwind difference scheme on a novel Shishkin--type mesh for a linear
convection--diffusion problem. {\em J.\ Comput.\ Appl.\ Math.} {\bf 110},
93--104 (1999).

\bibitem{rvpre}
R.~Vulanovi\'c, A priori meshes for singularly perturbed quasilinear
two--point boundary value problems.
{\em IMA J. Numer.\ Anal.} {\bf 21}, 349–-366 (2001).

\bibitem{rl}
H.-G.~Roos and T.~Lin{\ss},
Sufficient conditions for uniform convergence on layer adapted grids.
{\em Computing}\/ {\bf 64}, 27--45 (1999).

\bibitem{lrv}
T.~Lin{\ss}, H.-G.~Roos, and R.~Vulanovi\'c,
Uniform pointwise convergence on Shishkin--type meshes for quasilinear
convection--diffusion problems. {\em Preprint Math-NM-03-1999, T.~U.~Dresden}
(1999).

\bibitem{rvgs}
R.~Vulanovi\'c, Non--equidistant generalizations of the Gushchin--Shchen\-nikov
scheme. {\em Z. Angew.\ Math.\ Mech.} {\bf 67}, 625--632 (1987).

\bibitem{rveo}
R.~Vulanovi\'c, Some improvements of the non--equidistant Engquist--Osher
scheme. {\em Appl.\ Math.\ Comput.} {\bf 40}, 147--164 (1990).

\bibitem{sr}
M.~Stynes and H.-G.~Roos,
The midpoint upwind scheme. {\em Appl.\ Numer.\ Math.} {\bf 23}, 361--374
(1997).

\bibitem{tor2}
T.~Lin{\ss},
Uniform second order pointwise convergence of a finite difference discretization
for a quasilinear problem. {\em Preprint Math-NM-08-1999, T.~U.~Dresden} (1999).

\bibitem{lor}
J.~Lorenz,
Stability and monotonicity properties of stiff quasilinear boundary problems.
{\em Univ.\ u Novom Sadu Zb.\ Rad.\ Prirod.--Mat.\ Fak.\ Ser.\ Mat.}
{\bf 12}, 151--175 (1982).

\bibitem{kts} R.~B.~Kellogg and A.~Tsan, Analysis of some difference
approximations for a singular perturbation problem without turning points.
{\em Math.\ Comp.} {\bf 32}, 1025--1039 (1978).

\end{thebibliography}
\end{document}